\documentclass{article}
\usepackage{amsfonts}
\usepackage{amsthm}
\usepackage{amsmath}
\usepackage{amscd}
\usepackage{graphics}
\usepackage{graphicx}        % standard LaTeX graphics tool
\usepackage{array}
\usepackage{float}
\usepackage{authblk}
\title{Incorrect conclusions drawn for plausible looking diagrams}
\author{Protopapas Eleftherios}
\affil{National Technical University of Athens \\ School of Applied Mathematical and Physical Sciences \\Department of Mathematics\\
9 Heroon Polytechneiou st. \\GR-15 780 Athens, Greece\\
e-mail: lprotopapas@math.ntua.gr\\
orcid id:0000-0001-9569-695X} 
\date{}

\begin{document}

\maketitle

\begin{abstract}
In Mathematics is common to make a mistake and therefore a false conclusion arises.
%, which sometimes lead to conclusions that are not correct. 
In each case it is important to recognize the mistake in order to avoid a similar one in the future. Geometric figures provide decisive 
help in order to have a strict mathematical proof, but also can easily lead to wrong conclusions without a mathematical proof.
 In this paper, several incorrect conclusions drawn for plausible looking diagrams are presented, motivated by a well-known faulty model 
for measuring the length of a segment. Similar models that lead to a contradiction are developed and a model that leads to the correct result is derived.
The presented models prove the usefulness of paradoxes and can be implemented in a classroom in order to point out to students the significance of a strict mathematical proof 
as well as the construction of a correct mathematical model. The geometric nature of the problems provides the opportunity to use 
a dynamic geometric software.

\end{abstract}


{\bf{MOS subject classification:}} 97D50, 97G99.

\section{Introduction}
A principal component in mathematics education is problem solving \cite{vosk}. Students often face difficulties in mathematics problem 
solving struggling first for 
modeling the problem and later for solving it \cite{tamsub}. Productive struggle and learning mathematics are connected \cite{bic}. 
The struggle reduces whenever you can draw a diagram, in order to recognize the correct mathematical path for the mathematical proof. 
Unfortunately, the diagram could lead to wrong conclusions, 
due to mathematical mistakes or its misinterpretation.

Mathematical modeling is used to describe a phenomenon and it is crucial that its mathematical description is flawless. 
In the case where the model is not correct, the obtained results are often bizarre. In mathematics this is unacceptable and 
the mistake must be recognized and corrected. This is an important step, because the gained knowledge and experience will be used 
in the future in order to avoid such mistakes.

Mistakes occur from flawed reasoning, arithmetic errors or faulty logic. Two simple fault statements \cite{bunch} are  the following:
$$x=1 \Rightarrow x^2=x \Rightarrow x^2-1=x-1 \Rightarrow $$
$$\Rightarrow (x-1)(x+1)=x-1 \Rightarrow x+1=1 \Rightarrow x=0\Rightarrow 0=1$$  
and
$$-1=\left(\sqrt{-1} \right)^2=\sqrt{-1}\cdot \sqrt{-1}=\sqrt{(-1) \cdot (-1)}=\sqrt{1}=1.$$
Similar statements can be found in \cite{bunch}, \cite{clark}, \cite{klystap}, \cite{riddles}.

Sometimes the obtained faulty statements are called paradoxes, 
gaining great attention from many authors.
`A paradox is an item of information which contradicts with your current schemata,
thus creating a cognitive conflict disbalance' \cite{kondra}. 
There are two main reasons that lead to a paradox:
lack of essential logic or application of logic to a situation where it is not applicable. 
Kondratieva \cite{kondra} used paradoxes in order to understand Mathematics, 
discussing themes from Algebra, Set Theory, Analysis, etc. Achilles and a tortoise, dichotomy and the stadium are the most 
known Zeno's paradoxes \cite{donmez}. Cantor's and Russel's paradoxes are widely studied \cite{donmez}, \cite{kondra}. 
Paradoxes from Calculus and set theory are also studied 
\cite{kleiner}, while Bell used paradoxes to study oppositions in Mathematics \cite{bell}. 

It is useful to incorporate the study of paradoxes in mathematical classes \cite{kondra}. In our days their use is limited. Kleiner and 
Movshovitz-Hadar \cite{kleiner} suggested a role for paradoxes in the teaching and learning process, emphasizing that paradoxes create 
curiosity, motivation and a better learning environment. Along this line Kondratieva \cite{kondra} suggested that in order to correct 
students behavior it would help if the teacher comes up with an example of an error leading to a paradox.

In this paper we implemented Kondratieva's suggestion presenting several paradoxes that a diagram is employed. 
These paradoxes connect Geometry with 
other mathematical areas, providing significant didactic results \cite{kondber}. Since Geometry is employed the use of a dynamic geometry 
software can be implemented. Moreover explanations of the paradoxes are presented.
%In what follows motivated by a known fallacy model in Analysis for the calculation of a segment's length \cite{kondra}, 
Specifically, in section \ref{str1} we show strange results that occur from diagrams and in section \ref{str2} strange results 
that occur from faulty reasoning.
In section \ref{lathos} faulty models for calculating a segment's length are shown and a correct one
 is presented in section \ref{corr}. In section \ref{disc} discussion of these models is conducted, 
while in section \ref{conc} the conclusions are presented.

\section{Strange results that occur from diagrams}
\label{str1}
\subsection{Infinite length with finite area}
\begin{figure}[H]
  \centering
  \includegraphics[width=.99\linewidth]{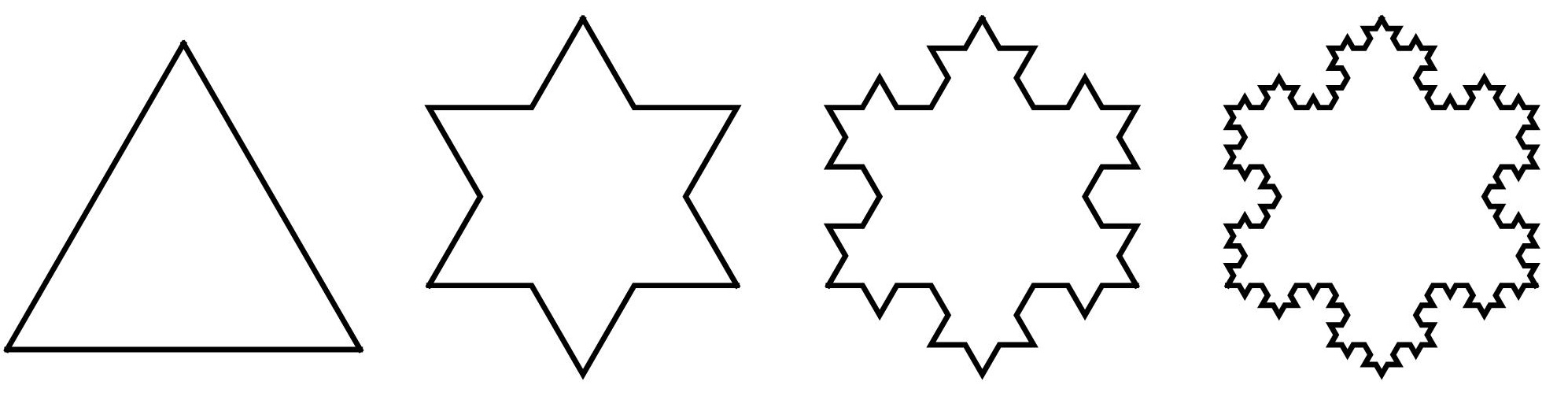}
   \caption{Koch snowflake.}
  \label{fig: koch1}
  \end{figure}

Let the length’s side of an equilateral triangle be $a>0.$ Each side is divided in three equal segments with side $\dfrac{a}{3}$ 
deleting the middle one. With side each one of the deleted ones we construct outside of the initial equilateral triangle three 
equilateral triangles with sides $\dfrac{a}{3}.$ Continuing this process indefinitely we draw the so-called 
Koch snowflake ({\bf{Figure \ref{fig: koch1}}}). The remarkable result of the Koch snowflake is that it has infinite length, 
but finite area \cite{koch}.

%Fractal geometry has many 
Regarding the perimeter, in the following table the values of interest are written.\vspace{1mm}
\\
\begin{tabular}{| m{1cm} | c | c | c |}
 \hline
Step & Equilateral(s) side length & Number of equal sides & Perimeter\\ 
\hline
 \\[-1em]
0 & $a$ & 3 & $3a$ \\ \hline 
\\[-1em]
 1 & $\dfrac{a}{3}$ & $3\cdot 4$ & $3\cdot 4 \dfrac{a}{3}$ \\ [0.5em] \hline 
\\[-1em]
2 & $\dfrac{a}{3^2}$ & $3\cdot 4^2$ & $3\cdot 4^2 \dfrac{a}{3^2}$ \\ [0.5em] \hline 
\\[-1em]
3 & $\dfrac{a}{3^3}$ & $3\cdot 4^3$ & $3\cdot 4^3 \dfrac{a}{3^3}$ \\ [0.5em] \hline 
\\[-1em]
\vdots &\vdots & \vdots & \vdots\\ \hline 
\\[-1em]
n & $\dfrac{a}{3^n}$ & $3\cdot 4^n$ & $3\cdot 4^n \dfrac{a}{3^n}$ \\ [0.5em] \hline 
\end{tabular}
\vspace{1mm} \\
From the previous table  the limit of the perimeter of the Koch snowflake when $n$ tends to infinity is indefinite, since
$$\lim_{n \rightarrow +\infty} \left( 3\cdot 4^n \dfrac{a}{3^n}\right)=3a\lim_{n \rightarrow +\infty} \left( \dfrac{4}{3}\right)^n=+\infty.$$
For completeness we have to prove the expression for the perimeter of the Koch snowflake using mathematical induction .

In order to calculate the area of the Koch snowflake, we create the following table containing all the necessary information.
\\
\begin{tabular}{| m{0.6cm} | c | c | c |}
 \hline
Step & Side's length & New triangles & Area \\ 
\hline
 \\[-1em]
0 & $a$ & 0 & $A_0=a^2\dfrac{\sqrt{3}}{4}$ \\ [0.5em]  \hline 
\\[-1em]
 1 & $\dfrac{a}{3}$ & 3 & $A_1=A_0+3\;\left(\dfrac{a}{3}\right)^2\;\dfrac{ \sqrt{3}}{4}$
\\ [0.5em] \hline \\[-1em]
2 & $\dfrac{a}{3^2}$ & $3\cdot 4$ & $A_2=A_1+3\cdot 4 \; \left(\dfrac{a}{3^2}\right)^2\;\dfrac{ \sqrt{3}}{4}$
\\ [0.5em] \hline \\[-1em]
3 & $\dfrac{a}{3^3}$ & $3\cdot 4^2$ & $A_3=A_2+3\cdot 4^2 \; \left(\dfrac{a}{3^3}\right)^2\;\dfrac{ \sqrt{3}}{4}$ 
\\ [0.5em] \hline \\[-1em]
\vdots &\vdots & \vdots & \vdots
\\ \hline \\[-1em]
n & $\dfrac{a}{3^n}$ & $3\cdot 4^n$ & $A_n=A_{n-1}+3\cdot 4^{n-1} \; \left(\dfrac{a}{3^{n}}\right)^2\;\dfrac{ \sqrt{3}}{4}$
\\ [0.5em] \hline 
\end{tabular}
\vspace{1mm} \\
It is obvious that
$$A_1=A_0\left(1+\frac{3}{9}\right),\;A_2=A_0\left(1+\frac{3}{9}+\frac{3\cdot 4}{9^2}\right),$$
$$A_3=A_0\left(1+\frac{3}{9}+\frac{3\cdot 4}{9^2}+\frac{3\cdot 4^2}{9^3}\right)$$
and therefore
$$A_n=A_0\left(1+\frac{3}{9}+\frac{3\cdot 4}{9^2}+...+\frac{3\cdot 4^{n-1}}{9^n}\right)=A_0\left(1+\sum_{k=1}^{n}\frac{3\cdot 4^{k-1}}{9^k}\right),\;n \in \mathbb{N}^*,$$
which can be proved using mathematical induction. 
This relation leads to the area of the Koch snowflake, which is
$$A=\lim_{n \rightarrow +\infty}A_n= A_0\left(1+\frac{\frac{3}{9}}{1-\frac{4}{9}}\right)=\frac{2\sqrt{3}a^2}{5},$$
since the series is geometric with ratio $\dfrac{4}{9}.$

	\subsection{Infinite area, but finite volume!}
\begin{figure}[H]
  \centering
  \includegraphics[width=.6\linewidth]{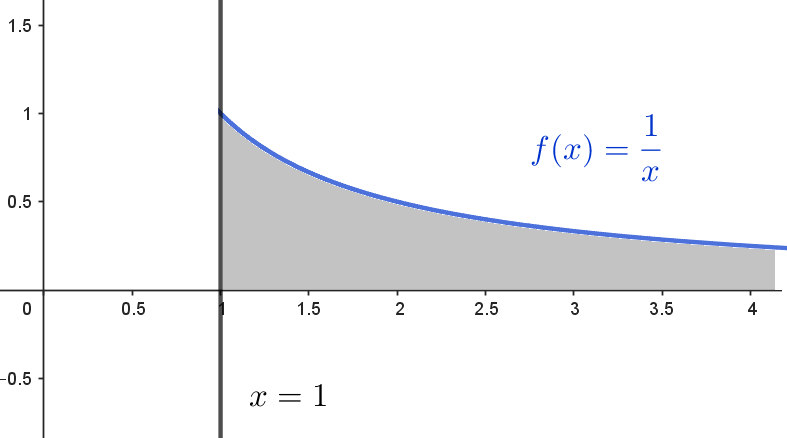}
   \caption{The graph of $f(x)=\dfrac{1}{x}.$}
  \label{fig: plo1_x}
  \end{figure}

Let $f$ be the function $f(x)=\dfrac{1}{x},\;x \in [1,+\infty).$ Since 
\begin{equation}
\int_1^{+\infty}f(x)dx=\lim_{A \rightarrow +\infty}\int_1^{A}\frac{1}{x}dx=\lim_{A \rightarrow +\infty}\left[ lnx \right]_1^{A}=\lim_{A \rightarrow +\infty}lnA=+\infty,
\end{equation}
the area between $C_f,$ axis $x'x$ and from the line $x=1$ up to $+\infty$ is indefinite ({\bf{Figure \ref{fig: plo1_x}}}).

\begin{figure}[H]
  \centering
  \includegraphics[width=.6\linewidth]{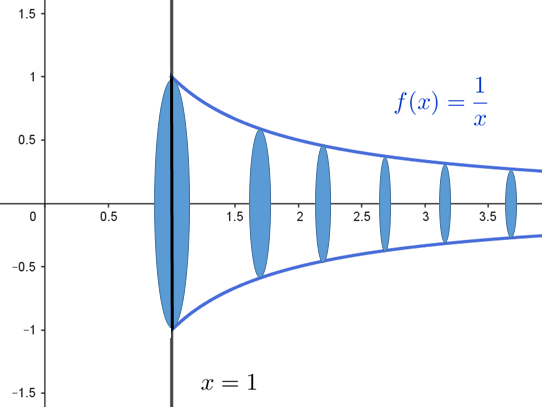}
   \caption{Rotating the graph of $f(x)=\dfrac{1}{x}$ along the x-axis.}
  \label{fig: plo1a_x}
  \end{figure}

It seems `logical' that rotating $C_f$ around x-axis the volume of the resulting solid of revolution would be also indefinite.
But, calculating the integral
$$
\int_1^{+\infty}f^2(x)dx= \lim_{A \rightarrow +\infty} \int_1^{A}\frac{1}{x^2}dx=\lim_{A \rightarrow +\infty}
\left[ -\frac{1}{x}\right]_1^{A}=\lim_{A \rightarrow +\infty}\left(-\frac{1}{A}+1 \right)=1,
$$ we conclude that rotating $C_f$ around x-axis the volume of the resulting solid of revolution is finite 
({\bf{Figure \ref{fig: plo1a_x}}}),  since
\begin{equation}
V=\int_1^{+\infty}\pi f^2(x)dx=\pi,
\end{equation}
which was unexpected.

	\section{Strange results that occur from faulty reasoning}
	\label{str2}
\subsection{The missing square}
\begin{figure}[H]
  \centering
  \includegraphics[width=.99\linewidth]{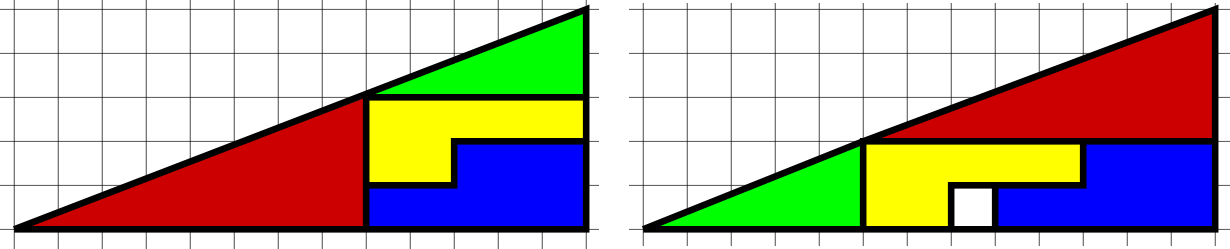}
   \caption{The missing square.}
  \label{fig: mis_squ}
  \end{figure}

One of the most famous paradoxes is the one described in {\bf{Figure \ref{fig: mis_squ}}}, where the two figures seem equal, 
but the areas are not, because of the missing square \cite{kondra}.
Strict geometrical analysis of the figures, reveal the faulty reasoning. Looking at both figures the assumption that there are two equal orthogonal triangles 
with area $\frac{13 \cdot 5}{2}\;cm^2=32.5\;cm^2$ seems logical, but this does not explain the missing square. Moreover, in both figures the colored part has area $32\;cm^2$ 
and therefore the problem remains unexplained, revealing though that the shapes are not orthogonal triangles. Calculating the slopes of 
the hypotenuses with respect to the horizontal sides of the red and the green orthogonal triangles,
we find that they are $\frac{3}{8},$ $\frac{2}{5},$ respectively, concluding that the hypotenuses of the triangles are not collinear. 
Changing the position of the red and the green triangle in the second
 figure the extra square appears, because the hypotenuses move a little bit up creating the space for our missing square. 
Drawing the figures with extreme caution and detail this can be easily verified.
\subsection{Bizarre areas}
\begin{figure}[H]
  \centering
  \includegraphics[width=.99\linewidth]{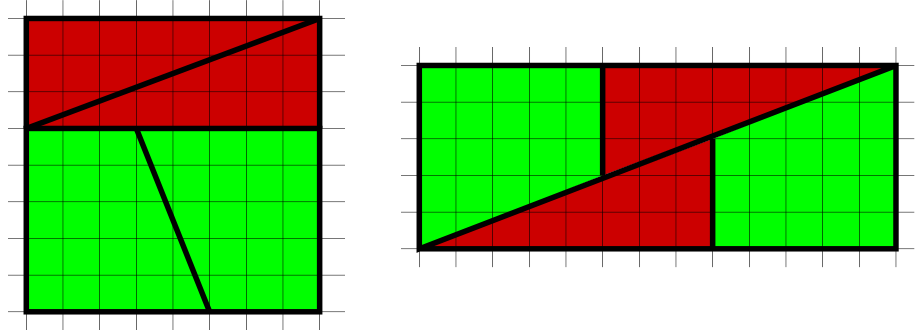}
   \caption{Bizarre areas.}
  \label{fig: areas1}
  \end{figure}

One more paradox is given in {\bf{Figure \ref{fig: areas1}}}, where both areas seem equal \cite{bunch}, but this means that $64=65!$
The mathematics reveal the unorthodox of the conclusion, which as in the previous subsection involves slopes of segments, 
justifying the fact that in the right figure the `diagonal' segments do not coincide, leaving a tiny area of $1\;cm^2$ uncolored, something that 
it is not visible in the figure. Once more cautious drawing will reveal the paradox.

\subsection{Aristotle's wheel paradox}
Every point of a rotating circle along a line without sliding belong to curve called a cycloid \cite{cajori}. Every point $(x,y)$ of a cycloid is defined via the relations
	\begin{equation}
	\begin{cases}
	x=R(\phi-\sin \phi)\\
	y=R(1-\cos \phi),
		\end{cases}
		\end{equation}
	where $\phi$ is the angle through which the rolling circle has rotated, the circle has radii $R$ and the point $O(0,0)$ belongs to the cycloid.\\

\begin{figure}[H]
  \centering
  \includegraphics[width=.9\linewidth]{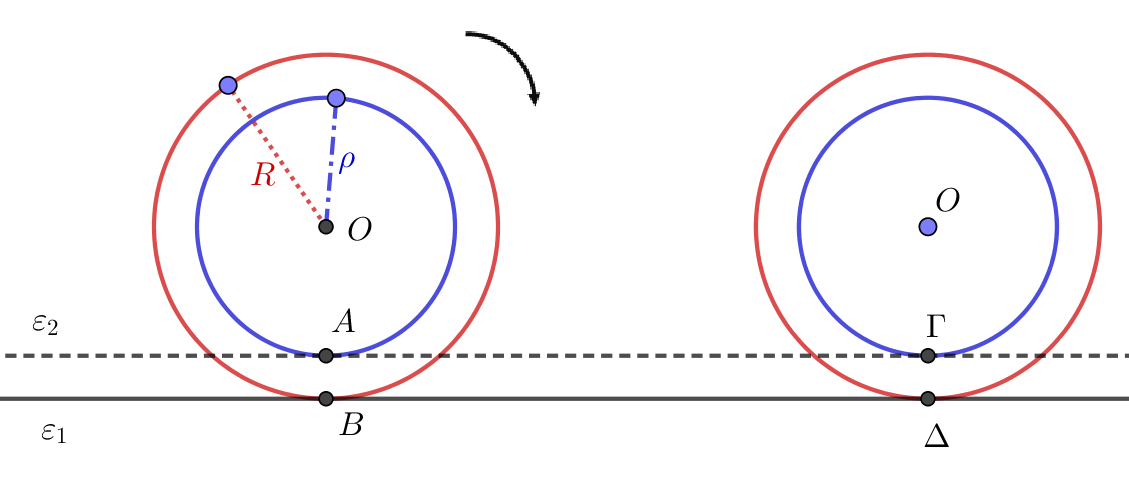}
   \caption{Rotating circles.}
  \label{fig: circle}
  \end{figure}

Consider the following problem. Two circles with common center $O$ and different radii $R > \rho > 0$ are attached to each other. 
The parallel lines $\varepsilon_1,$ $\varepsilon_2$ 
are tangents of the circles $(O,R),$ $(O,\rho)$ in the points $B,$ $A,$ respectively \cite{klystap}.
These circles rotate clockwise along their tangent lines $\varepsilon_1,$ $\varepsilon_2$ without sliding and complete a complete rotation 
reaching at the points
 $\Delta,$ $\Gamma,$ respectively ({\bf{Figure \ref{fig: circle}}}). 
Since $B \Delta = A\Gamma,$ we conclude that 
$$2 \pi R=2 \pi \rho \Leftrightarrow R=\rho!$$
	
	It is an extraordinary result, but Mathematics reveal the faulty reasoning of this problem. 
	 		If the circles make a full turn without sliding, then $R=\rho,$ which is a contradiction. Therefore there is a sliding circle. 
			It is easy to realize that the inner circle is the sliding one. As a matter of fact the smaller the inner circle is the more it slides. 
			This result can be verified using a dynamic geometric software, in which we can see that the inner circle's point do not belong to a 
			cycloid.
		
	%Therefore, if $(x_A,y_A),$ $(x_B,y_B),$ 
	%are the coordinates of $A,$ $B,$ respectively, we reach at
	%\begin{equation}
	%\begin{cases}
	%x_A=R\left(\phi+\frac{\pi}{2}-\cos \phi\right)\\
	%y_A=R\sin \phi,
		%\end{cases}
	%\end{equation}
	%\begin{equation}
	%\begin{cases}
	%x_B=\rho \left(\phi+\frac{\pi}{2}-\cos \phi\right)\\
	%y_B=\rho \sin \phi,
		%\end{cases}
	%\end{equation}
	%where $\phi \in \left[-\dfrac{\pi}{2},\dfrac{3\pi}{2}\right].$
	\section{Measuring the length of a segment}
\label{lathos}
In this section four models are presented in order to derive the length of a segment, which are self-similar 
denoting a fractal structure with fractal dimension equal to 1 \cite{mand2}. The first one is a well known model, while all the others 
are developed in order to create new paradoxes as well as a correct model.
 In each faulty model a similar mistake is embedded 
and an extraordinary result is derived. 
%In subsection \ref{pi} a well know model is presented, in subsection \ref{riza} a similar 
%model is shown, while in subsection \ref{lamda} a generalized model is studied. In subsection \ref{diafora} two more faulty models are presented.
%we discuss those three models providing answers of the obtained unusual results.

\subsection{Is $2$ equals to $\pi?$}
\label{pi}
\begin{figure}[H]
\begin{minipage}{.49\textwidth}
  \centering
  \includegraphics[width=.99\linewidth]{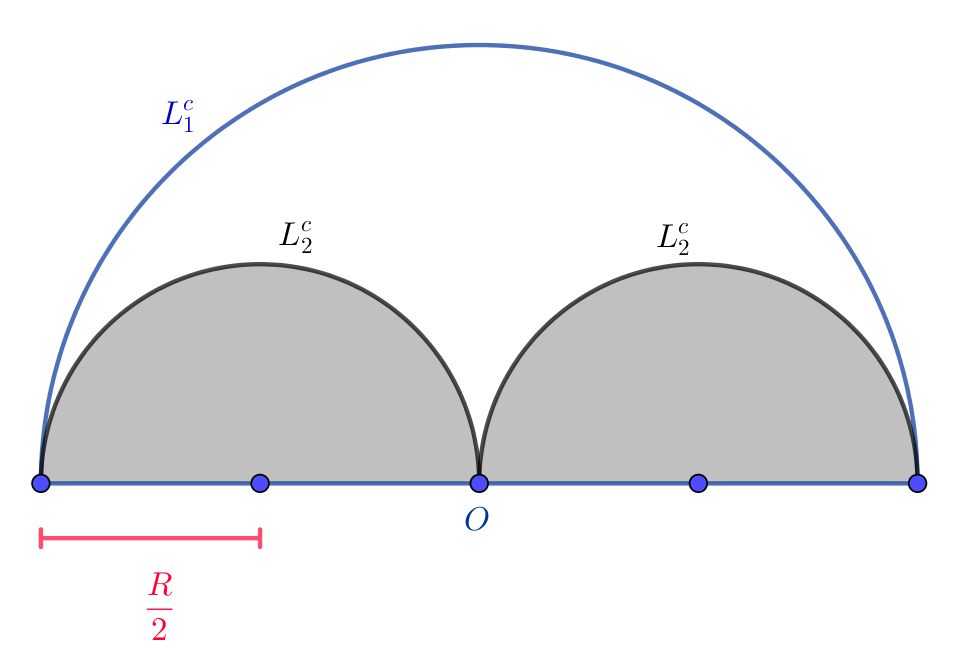}
   \caption{The second approach with semicirles.}
  \label{fig: K2}
  \end{minipage}
\begin{minipage}{.02\textwidth}
   \end{minipage}
\begin{minipage}{.49\textwidth}
 \centering
  \includegraphics[width=.99\linewidth]{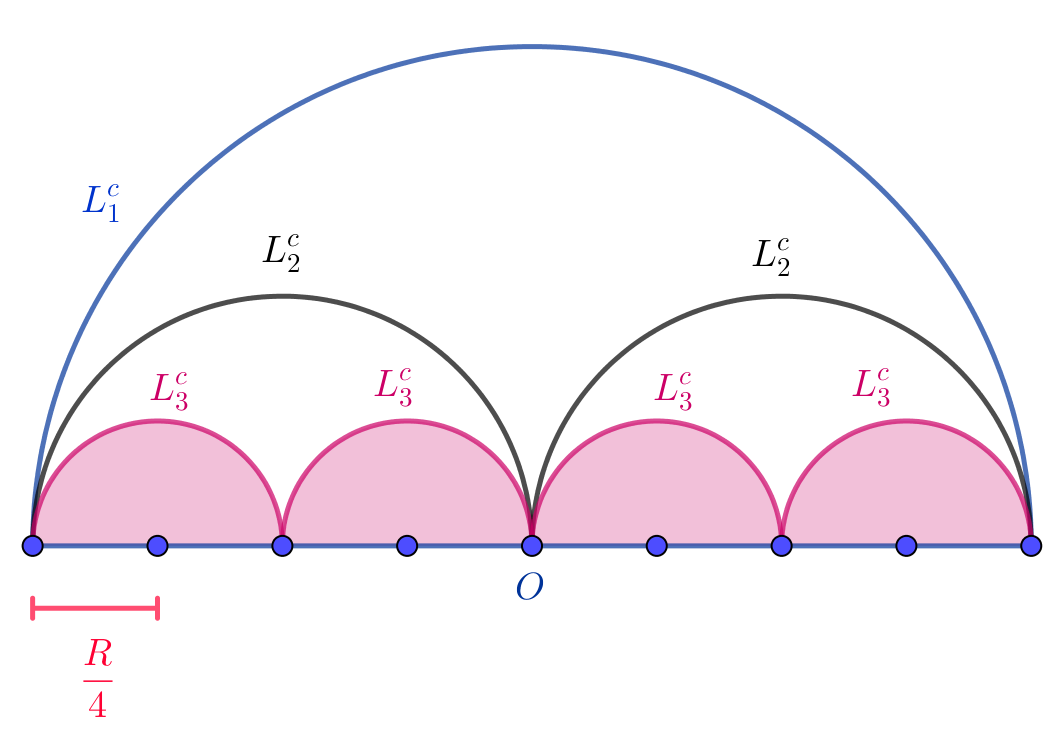}
   \caption{The third approach with semicirles.}
  \label{fig: K3}
  \end{minipage}
\end{figure}

Let $O$ be the center of a semicircle with radius $R>0$ 
%({\bf{Figure \ref{fig: K1}}}) 
and perimeter 
\begin{equation}
\label{eqn: L1}
L_1^c=\pi R.
\end{equation}
%\begin{figure}[H]
 %\centering
 %\includegraphics[width=.6\linewidth]{K1}
%\caption{The first approach with a semicirle.}
 %\label{fig: K1}
%\end{figure}
With diameter $R,$ two new semicircles are drawn ({\bf{Figure \ref{fig: K2}}}), where each one has perimeter
\begin{equation}
\label{eqn: L2}
L_2^c=\pi \frac{R}{2}.
\end{equation}
Continuing this procedure four semicircles are drawn ({\bf{Figure \ref{fig: K3}}}) and 
after  $n \in \mathbb{N}^*$ repetitions, $2^{n-1}$ semicircles are drawn with perimeter
\begin{equation}
\label{eqn: Ln}
L_n^c=\pi \frac{R}{2^{n-1}},\;n \in \mathbb{N}^*
\end{equation}
and therefore all the semi-circumferences add up to 
\begin{equation}
\label{eqn: Sn}
S_n^c=2^{n-1}L_n^c=2^{n-1}\pi \frac{R}{2^{n-1}}=\pi R.
\end{equation}
When $n$ tends to infinity, the sum, $S_n^c,$ `tends' to the diameter $2R$ and since
\begin{equation}
\label{eqn: S1}
\lim_{n \rightarrow +\infty}S_n^c=\lim_{n \rightarrow +\infty}\pi R=\pi R,
\end{equation}
it stands
\begin{equation}
\label{eqn: S2}
2R=\pi R \Leftrightarrow \pi=2!
\end{equation}

\subsection{Is $2$ equals to $\sqrt{2}?$}
\label{riza}
\begin{figure}[H]
\begin{minipage}{.49\textwidth}
  \centering
  \includegraphics[width=.99\linewidth]{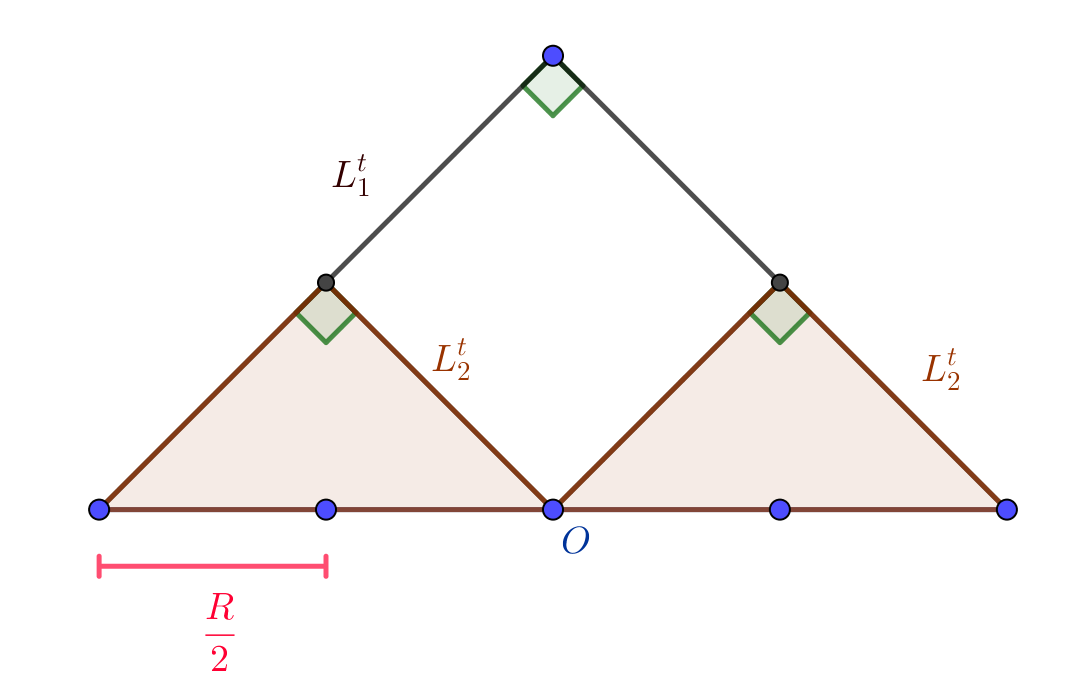}
   \caption{The second approach with isosceles orthogonal triangles.}
  \label{fig: KK2}
  \end{minipage}
\begin{minipage}{.02\textwidth}
   \end{minipage}
\begin{minipage}{.49\textwidth}
 \centering
  \includegraphics[width=.99\linewidth]{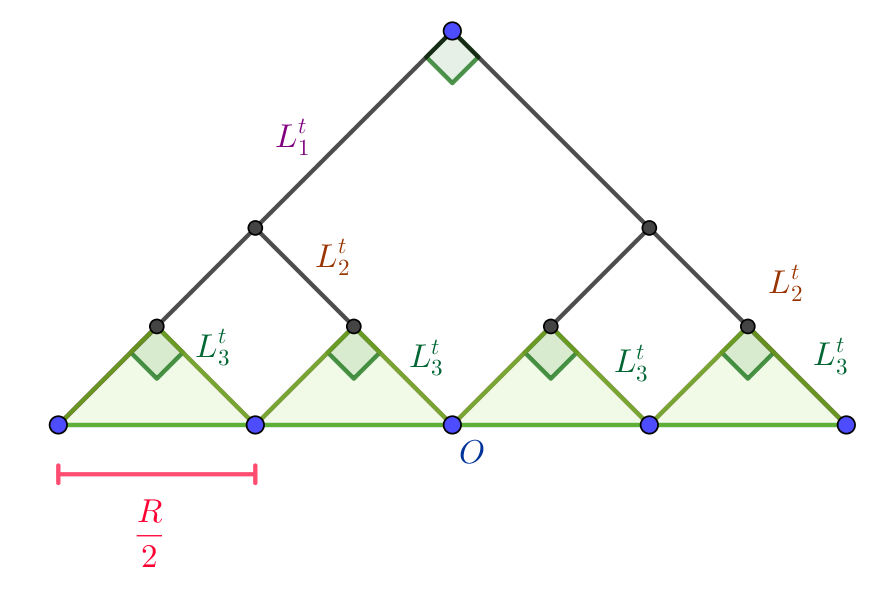}
   \caption{The third approach with isosceles orthogonal triangles.}
  \label{fig: KK3}
  \end{minipage}
\end{figure}

Let an isosceles orthogonal triangle with hypotenuse's length $2R,$ $R>0,$ 
%({\bf{Figure \ref{fig: KK1}}}) 
and sum of the equals sides 
\begin{equation}
\label{eqn: LL1}
L_1^t=2\sqrt{2}R.
\end{equation}
With hypotenuse equal to $R,$ two isosceles orthogonal triangles are drawn ({\bf{Figure \ref{fig: KK2}}}), 
where each one has a sum of the equal sides
\begin{equation}
\label{eqn: LL2}
L_2^t=\sqrt{2}R.
\end{equation}

Continuing the same procedure $n \in \mathbb{N}^*$ times (see for instance {\bf{Figure \ref{fig: KK3}}}), the sum of the equal sides 
of each of the 
$2^{n-1}$ isosceles triangles is
\begin{equation}
\label{eqn: LLn}
L_n^t=\sqrt{2}\frac{R}{2^{n-2}},
\end{equation}
while the total sum of the equal sides of the $2^{n-1}$ isosceles triangles is
\begin{equation}
\label{eqn: SS1}
S_n^t=2^{n-1}L_n^t=2^{n-1}\sqrt{2}\frac{R}{2^{n-2}}=2\sqrt{2}R,
\end{equation}
so
\begin{equation}
\label{eqn: SS1}
\lim_{n \rightarrow +\infty}S_n^t=\lim_{n \rightarrow +\infty}2\sqrt{2} R=2\sqrt{2} R,
\end{equation}
which is the diameter, so
\begin{equation}
\label{eqn: S2}
2R=2\sqrt{2} R \Leftrightarrow \sqrt{2}=2!
\end{equation}

\subsection{Is $2$ equals to a number in $(2,2\sqrt{2}]?$}
\label{lamda} 

In the previous model the ratio of the vertical sides of each isosceles triangle is $1.$ 
In this section the assumption is that this ratio 
is $\lambda>0,$ so when the hypotenuse is $2R,$ $R>0,$ the vertical segments are 
\begin{equation}
\label{eqn: seg}
\frac{2R}{\sqrt{\lambda^2+1}},\;\;\frac{2 \lambda R}{\sqrt{\lambda^2+1}}.
\end{equation}
This assumption leads to a model in which similar orthogonal triangles with similarity ratio equals to $\lambda$ are used 
between steps and in every step of the construction. 

\begin{figure}[H]
\begin{minipage}{.49\textwidth}
  \centering
  \includegraphics[width=.99\linewidth]{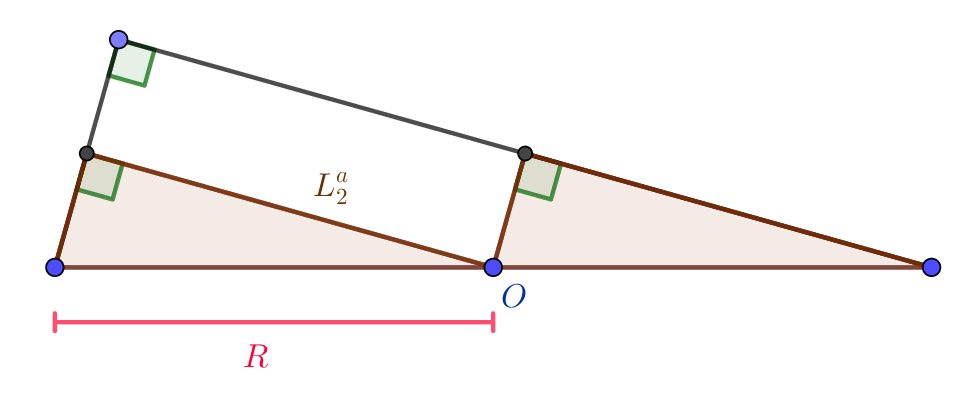}
   \caption{The second approach with similar orthogonal triangles.}
  \label{fig: or82}
  \end{minipage}
\begin{minipage}{.02\textwidth}
   \end{minipage}
\begin{minipage}{.49\textwidth}
 \centering
  \includegraphics[width=.99\linewidth]{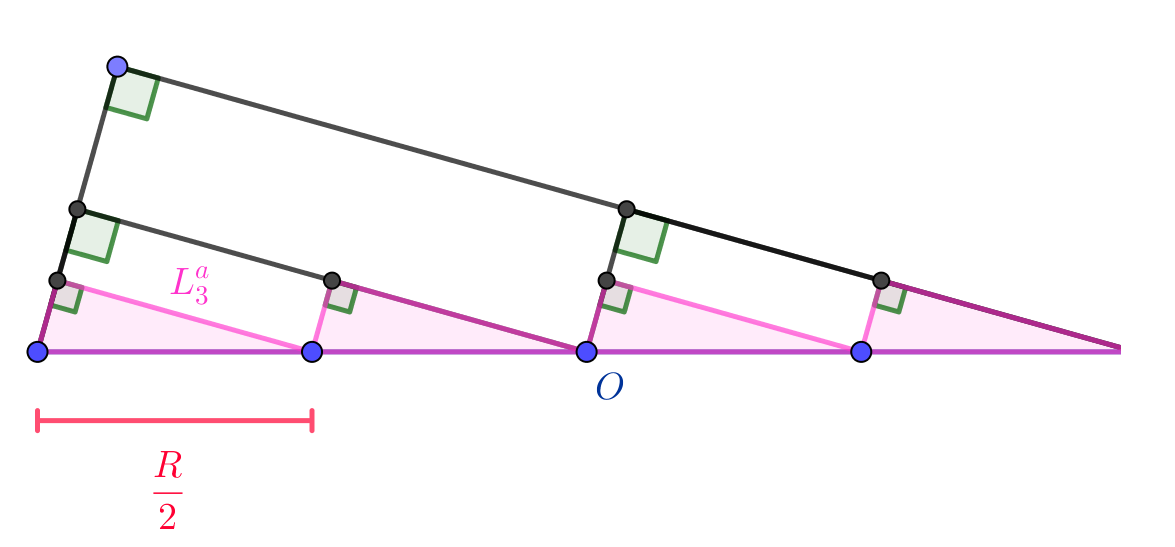}
   \caption{The third approach with similar orthogonal triangles.}
  \label{fig: or83}
  \end{minipage}
\end{figure}

Following the procedure (after $n \in \mathbb{N}^* $ repetitions) that is described in section \ref{pi}
(see consequently 
%{\bf{Figure \ref{fig: or81}}}, 
{\bf{Figure \ref{fig: or82}}}, {\bf{Figure \ref{fig: or83}}}), 
the total sum of the vertical segments of the 
$2^{n-1}$ triangles is
\begin{equation}
\label{eqn: seg1}
S_n^{\lambda}=\frac{2(1+\lambda)R}{\sqrt{\lambda^2+1}}
\end{equation}
and therefore
\begin{equation}
\label{eqn: seg2}
\frac{2(1+\lambda)R}{\sqrt{\lambda^2+1}}=2R \Leftrightarrow 2=\frac{2(1+\lambda)}{\sqrt{\lambda^2+1}}.
\end{equation}
Function $$f:D_f=(0,+\infty) \rightarrow \mathbb{R}$$ with
\begin{equation}
\label{eqn: func1}
f(x)=\frac{2(1+x)}{\sqrt{x^2+1}},\;x>0,
\end{equation}
is continous in $(0,+\infty),$ strictly increasing in the interval $(0,1],$ strictly decreasing in the interval $[1,+\infty)$
and 
 \begin{equation}
\label{eqn: func2}
f\big(D_f\big)=\left(2,2\sqrt{2}\;\right].
\end{equation}

Therefore from (\ref{eqn: seg2}) every number in the interval $\left(2,2\sqrt{2}\;\right]$ must be equal to $2!$

\subsection{What is really equals to 2?}
\label{diafora}
\begin{figure}[H]
  \centering
  \includegraphics[width=.5\linewidth]{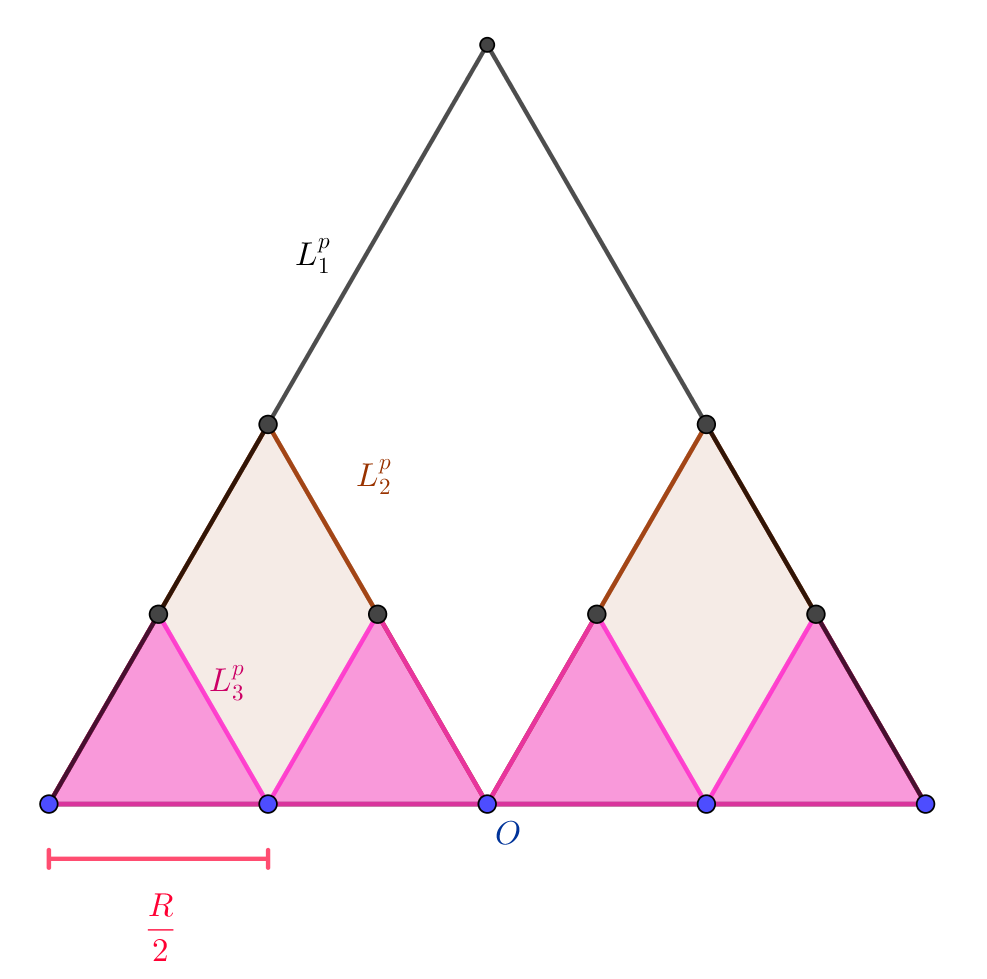}
   \caption{The third approach with equilateral triangles.}
  \label{fig: isopleuro}
  \end{figure}

The procedure that is described previously can be applied with equilateral triangles 
({\bf{Figure \ref{fig: isopleuro}}}).

It is easy to verify that in this case
\begin{equation}
2=4!
\end{equation}

\section{Constructing a correct sequence}
\label{corr}
In order to construct a sequence that converges to the correct value, triangles that are not similar in every step of the procedure are used. 
For convenience isosceles triangles are used, with decreasing equal angles.

Let a segment with length $2R>0$ and an isosceles triangle where the equal angles are $\omega,$ the side which is unequal 
to the other two is $2R$ 
%({\bf{Figure \ref{fig: iso1}}}) 
and therefore the sum of its equal sides is
\begin{equation}
\label{eqn: SSSS1}
S_1^k=L_1^k=\frac{2R}{cos(\omega)}.
\end{equation}

\begin{figure}[H]
\begin{minipage}{.49\textwidth}
  \centering
  \includegraphics[width=.99\linewidth]{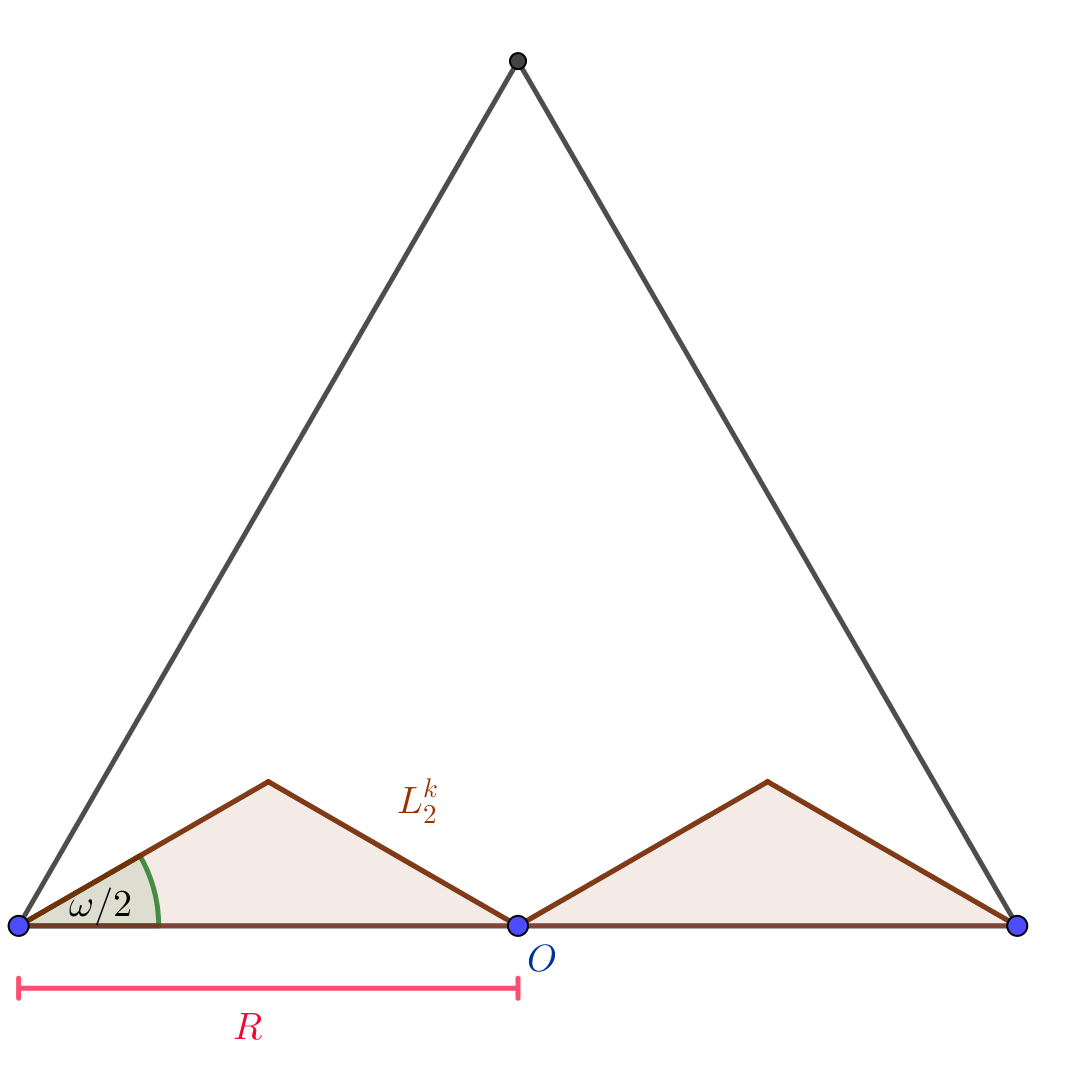}
   \caption{The second approach with non similar isosceles triangles.}
  \label{fig: iso2}
  \end{minipage}
\begin{minipage}{.02\textwidth}
   \end{minipage}
\begin{minipage}{.49\textwidth}
 \centering
  \includegraphics[width=.99\linewidth]{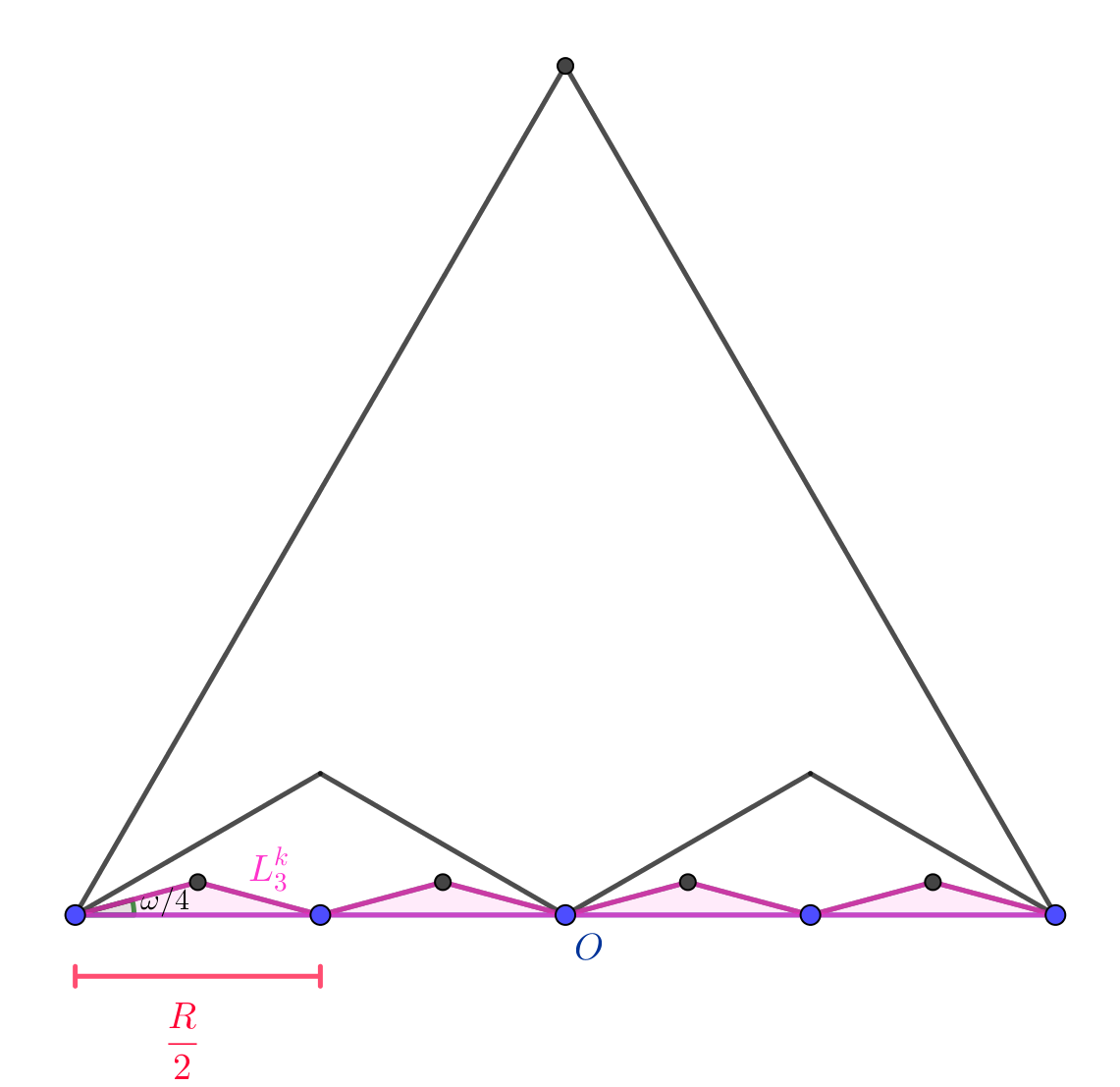}
   \caption{The third approach with non similar isosceles triangles.}
  \label{fig: iso3}
  \end{minipage}
\end{figure}

Bisecting angle $\omega$ and taking the base of the isosceles triangle equals to $R$ ({\bf{Figure \ref{fig: iso2}}}), 
the total sum of the equal sides of each triangle is
\begin{equation}
\label{eqn: SSSS2}
S_2^k=2L_2^k=\frac{2R}{cos\left(\frac{\omega}{2}\right)}.
\end{equation}

Therefore in the $n-th$ repetition $(n \in \mathbb{N}^*)$ of this procedure (see {\bf{Figure \ref{fig: iso3}}})
the total sum of the $2^{n-1}$ triangles' equal sides  is
\begin{equation}
\label{eqn: SSSSn}
S_n^k=2^{n-1}L_n^k=\frac{2R}{cos\left( \frac{\omega}{2^{n-1}} \right )},
\end{equation}
so
\begin{equation}
\label{eqn: SSSnlim}
\lim_{n \rightarrow +\infty}S_n^c=\lim_{n \rightarrow +\infty}\frac{2R}{cos\left(\frac{\omega}{2^{n-1}}\right)}=2R.
\end{equation}

\section{Discussion}
\label{disc}
It is worth wondering why the models from section \ref{lathos} came up with these results. 
It is obvious that sequences, $S_n,$ are constant ($\pi R,$ $2\sqrt{2}R,$ $f(\lambda)R,$ and $4,$ respectively), 
and therefore they converge in their constant value, which by definition is not equal to $2R.$ 
This is due to the fact that in each step the 
figures are similar to each other, creating fractal images. For instance, looking at the models that triangles are used, as the segment decreases two times, 
the triangle’s sides decrease two times.
The fractal nature of the curves involved 
denote that if they are observed under a microscope, they still would have the same shape. On the contrary with the fractal 
Koch snowflake \cite{koch} which has finite area and infinite perimeter at the same time, the derived structures have finite length and area.

In the case of the semicircles, since the semi-circumference (section \ref{pi}) is greater than the diameter in the $n-th$ step it stands 
for the semi-circumference that
$$\pi \frac{R}{2^{n-1}}>\frac{R}{2^{n-1}},$$
while in the triangles' models (section \ref{riza})
 from the triangle's inequality is obtained
$$2\sqrt{2}\frac{R}{2^{n-1}} > 2\frac{R}{2^{n-1}},\;f(\lambda)\frac{R}{2^{n-1}} > 2\frac{R}{2^{n-1}},\;2 \cdot 4\frac{R}{2^{n-1}}>2 \frac{R}{2^{n-1}},$$ 
respectively, indicating the faulty models. 

On the contrary, the triangle's inequality in the model that converges in the 
correct value (section \ref{corr}) is given by 
$$2\frac{\frac{R}{2^{n-1}}}{cos\frac{\omega}{2^{n-1}}}>\frac{R}{2^{n-1}},$$
and 
$$\lim_{n \rightarrow +\infty}\left(2\frac{\frac{R}{2^{n-1}}}{cos\frac{\omega}{2^{n-1}}}-\frac{R}{2^{n-1}}\right)=0,$$
verifying the convergence of the derived model.

\section{Conclusions}
\label{conc}
In this paper we presented several mathematical fallacies. These fallacies occur taking into account only the diagrams, creating
wrong conclusions. When a strict mathematical proof is employed the fallacies uncover, revealing the correct answer. 
Sometimes drawing a detailed diagram  using a dynamic 
 geometric software the mistake appears. It is usual to use diagrams 
in order to compose a strict mathematical proof. Students have to recognize the necessity of the strict mathematical procedure, 
which in many cases can be 
easier developed using a diagram.

The fallacies that are presented  for measuring the length of a segment prove the value of the correct definition of a mathematical model. 
From these paradoxes it can be derived that something may seem visually correct in geometry, but this does not mean that 
is correct or can be proved mathematically! Areas and perimeters behave very differently when 
limits are concerned! In the corresponding presented paradoxes the length of a limit curve is not the same as the limit of the length of the curve.

Careful designing and proper investigation are necessary in order to avoid mistakes that are crucial. Mistakes generate curiosity, 
increase motivation, create a learning environment  and show that faulty logic is common in Mathematics, but can be useful for 
Mathematical progress.

Teachers must implement paradoxes in their teaching, because the strange result increases student's interest, leading to new learning paths 
in order to gain knowledge. This will improve the problem-solving ability of the students. The presented paradoxes establish bridges between 
Geometry and other mathematical areas emphasizing in the necessity of a strict mathematical proof in order to verify or not what it seems logical.

\section*{Acknowledgment}
The author would like to thank M.Sc. mathematicians Kyriazis Christos and Sabani Maria, because the idea and the motivation for this paper came 
from fruitful discussions with them.

\end{document}